\documentclass{amsart}
\usepackage{amssymb}
\usepackage{graphicx}

\newtheorem{thm}{Theorem}
\newtheorem{cor}{Corollary}
\newtheorem{lem}{Lemma}

\newtheorem{rem}{Remark}

\newtheorem{conj}{Conjecture}
\theoremstyle{definition}
\newtheorem{example}[equation]{Example}
\newtheorem{prob}[equation]{Problem}

\newcommand{\ID}{{\mathbb D}}

\newcommand{\D}{{\mathbb D}}

\def\be{\begin{equation}}
\def\ee{\end{equation}}

\newcommand{\bee}{\begin{enumerate}}
\newcommand{\eee}{\end{enumerate}}

\newcommand{\blem}{\begin{lem}}
\newcommand{\elem}{\end{lem}}
\newcommand{\bthm}{\begin{thm}}
\newcommand{\ethm}{\end{thm}}
\newcommand{\bcor}{\begin{cor}}
\newcommand{\ecor}{\end{cor}}
\newcommand{\beg}{\begin{example}}
\newcommand{\eeg}{\end{example}}
\newcommand{\begs}{\begin{examples}}
\newcommand{\eegs}{\end{examples}}
\newcommand{\bdefe}{\begin{defin}}
\newcommand{\edefe}{\end{defin}}
\newcommand{\bprob}{\begin{prob}}
\newcommand{\eprob}{\end{prob}}
\newcommand{\bei}{\begin{itemize}}
\newcommand{\eei}{\end{itemize}}

\newcommand{\bcon}{\begin{conj}}
\newcommand{\econ}{\end{conj}}
\newcommand{\bcons}{\begin{conjs}}
\newcommand{\econs}{\end{conjs}}
\newcommand{\bprop}{\begin{propo}}
\newcommand{\eprop}{\end{propo}}
\newcommand{\br}{\begin{rem}}
\newcommand{\er}{\end{rem}}
\newcommand{\brs}{\begin{rems}}
\newcommand{\ers}{\end{rems}}
\newcommand{\bo}{\begin{obser}}
\newcommand{\eo}{\end{obser}}
\newcommand{\bos}{\begin{obsers}}
\newcommand{\eos}{\end{obsers}}
\newcommand{\bpf}{\begin{pf}}
\newcommand{\epf}{\end{pf}}
\newcommand{\ba}{\begin{array}}
\newcommand{\ea}{\end{array}}
\newcommand{\beq}{\begin{eqnarray}}
\newcommand{\beqq}{\begin{eqnarray*}}
\newcommand{\eeq}{\end{eqnarray}}
\newcommand{\eeqq}{\end{eqnarray*}}

\begin{document}

\title[  Hankel determinant for inverse of univalent functions]{ Hankel determinant of type $\boldsymbol{H_{2}(3)}$ for inverse functions of some classes of univalent
functions with missing second coefficient}

\author[M. Obradovi\'{c}]{Milutin Obradovi\'{c}}
\address{Department of Mathematics,
Faculty of Civil Engineering, University of Belgrade,
Bulevar Kralja Aleksandra 73, 11000, Belgrade, Serbia.}
\email{obrad@grf.bg.ac.rs}

\author[N. Tuneski]{Nikola Tuneski}
\address{Department of Mathematics and Informatics, Faculty of Mechanical Engineering, Ss. Cyril and
Methodius
University in Skopje, Karpo\v{s} II b.b., 1000 Skopje, Republic of North Macedonia.}
\email{nikola.tuneski@mf.ukim.edu.mk}

\subjclass{30C45, 30C55}

\keywords{Hankel determinant, inverse function, classes of univalent functions}

\begin{abstract}
In this paper we determine the upper bounds of $|H_{2}(3)|$ for the inverse functions of functions of some classes of univalent functions, where $H_{2}(3)(f)=a_{3}a_{5}-a_{4}^{2}$ is the Hankel determinant of a special type.
\end{abstract}

\maketitle

\section{Introduction and preliminaries}

Let ${\mathcal A}$ be the class containing functions that are analytic in the unit disk $\D:= \{ |z| < 1 \}$
and  are normalized such that
$f(0)=0= f'(0)-1$, i.e.,
\begin{equation}\label{eq-1}
f(z)=z+a_2z^2+a_3z^3+\cdots.
\end{equation}
By ${\mathcal S}$ we denote the class of functions from ${\mathcal A}$ which are univalent in $\D$.

Also, we need the classes of functions of bounded turning,  of convex functions, of starlike  functions, and of functions starlike with
respect to symmetric points, subclasses of ${\mathcal S}$, defined respectively in the following way
\[
\begin{split}
\mathcal{R}&=\left[f\in\mathcal{A}: {\rm Re}f'(z)>0,\, z\in\ID \right],\\
\mathcal{C}&=\left[f\in\mathcal{A}: {\rm Re}\left[1+\frac{zf''(z)}{f'(z)}\right]>0,\, z\in\ID \right],\\
\mathcal{S}^{\star}&=\left[f\in\mathcal{A}: {\rm Re}\frac{zf'(z)}{f(z)}>0,\, z\in\ID \right],\\
\mathcal{S}_{s}^{\star}&=\left[f\in\mathcal{A}: {\rm Re}\frac{2zf'(z)}{f(z)-f(-z)}>0,\, z\in\ID\right].
\end{split}
\]

\medskip

In his paper \cite{zaprawa} Zaprawa considered the Hankel determinant of the type
$$H_{2}(3)(f)=a_{3}a_{5}-a_{4}^{2},$$
defined for the coefficients of the function $f$ given by \eqref{eq-1}.
The author treated bounds of $|H_{2}(3)(f)|$ for the classes $\mathcal{R},\mathcal{C},\mathcal{S}^{\star}$
and gave sharp results in the case $a_{2}=0$. He also investigated the general case of these classes. In the
same paper it is proved that
$$ \max\left\{|H_{2}(3)(f)| : f\in \mathcal{S}\right\} >1.$$

\medskip

The object of current paper is to obtained the bounds of the modulus of the Hankel determinant  $H_{2}(3)(f^{-1})$ of coefficients of the inverse of function from the classes $\mathcal{R},$ $\mathcal{C},$ $\mathcal{S}^{\star}$ and $\mathcal{S}_{s}^{\star}$, defined before, as well as for the class
$\mathcal{S}$.
In all cases we suppose that function $f$ is missing its second coefficient, i.e., $a_{2}=0$.

\medskip

Namely, for every univalent function in $\D$ exists inverse at least on the disk with radius 1/4 (due to
the famous Koebe's 1/4 theorem). If the inverse has an expansion
\begin{equation}\label{eq-2}
f^{-1}(w) = w+A_2w^2+A_3w^3+\cdots,
\end{equation}
then, by using the identity $f(f^{-1}(w))=w$, from \eqref{eq-1} and \eqref{eq-2} we receive
\[
\begin{split}
A_{2}&=-a_{2}, \\
A _{3}&=-a_{3}+2a_{2}^{2} , \\
A_{4}&= -a_{4}+5a_{2}a_{3}-5a_{2}^{3},\\
A_{5}&= -a_{5}+6a_{2}a_{4}-21a_{2}^{2}a_3+3a_3^2+14a_2^4.
\end{split}
\]
Especially, when $a_{2}=0$, we have
\[
A_{2}=0, \quad A _{3}=-a_{3},\quad A_{4}= -a_{4},\quad A_{5}= -a_{5}+3a_3^2.
\]
So, in this case,
\begin{equation}\label{eq-5}
H_{2}(3)(f^{-1})=A _{3}A _{5}-A _{4}^{2}=a_{3}a_{5}-a_{4}^{2}-3a_3^3,
\ee
i.e.,
\begin{equation}\label{eq-6}
H_{2}(3)(f^{-1})=H_{2}(3)(f)-3a_3^3.
\ee

\medskip

For our further consideration we  need the next lemma given by Carlson \cite{carlson}.
\begin{lem}\label{lem-carl}
Let
\begin{equation}\label{eq-7}
\omega(z)=c_{1}z+c_{2}z^{2}+\cdots
\end{equation}
be a Schwartz function, i.e., a function analytic in $\D$, $\omega(0) = 0$ and $|\omega(z)| < 1$. Then
 $$|c_1|\le 1,\quad |c_2|\le1-|c_1|^2,\quad  |c_3|\le 1-|c_1|^2-\frac{|c_2|^2}{1+|c_1|},
 \quad \mbox{and}\quad |c_4|\le 1-|c_1|^2-|c_2|^2.$$
\end{lem}

\medskip

\section{Main results}

\bthm\label{22-th-1}
Let $f\in\mathcal{A}$ is given by \eqref{eq-1} and let $a_{2}=0$. Then
\begin{itemize}
\item[(a)]  $|H_{3}(2)(f^{-1})|\leq\frac{28}{45} $ if $f\in\mathcal{R}$;
\item[(b)]  $|H_{3}(2)(f^{-1})|\leq\frac{2}{45}  $ if $f\in\mathcal{C}$;
\item[(c)]  $|H_{3}(2)(f^{-1})|\leq 2$ if $f\in\mathcal{S}^{\star}$;
\item[(d)]  $|H_{3}(2)(f^{-1})|\leq 2$ if $f\in\mathcal{S}^{\star}_{s}$.
\end{itemize}
All these results are sharp.
\ethm

\begin{proof}$ $
\begin{itemize}
\item[(a)]
Since $f\in \mathcal{R}$ is equivalent to
$$f'(z)=\frac{1+\omega(z)}{1-\omega(z)},$$
for certain Schwartz function $\omega$, we receive that
\be\label{eq-8}
f'(z)=1+2\omega(z)+2\omega^{2}(z)+\cdots.
\ee
Using the notations for $f$ and $\omega$ given by \eqref{eq-1} and \eqref{eq-7}, and equating the coefficients in \eqref{eq-8}, we receive
\be\label{eq-9}
 \left\{\begin{split}
  a_2=& c_1, \\
  a_3=& \frac{2}{3}(c_2+c_1^2), \\
  a_4=& \frac{1}{2}(c_3+2c_1c_2+c_1^3), \\
  a_5=& \frac{2}{5}(c_4+2c_1c_3+ 3c_1^2c_2+c_2^2+c_1^4).
\end{split}
\right.
\ee
Since $a_{2}=0$, by \eqref{eq-9} we have $c_{1}=0$, and the appropriate coefficients have the next form:
\be\label{eq-10}
a_3= \frac{2}{3}c_2,\quad a_4= \frac{1}{2}c_3,\quad  a_5= \frac{2}{5}(c_4+c_2^2).
\ee
Now, from  \eqref{eq-5} and \eqref{eq-10}, after simple computation, we obtain
$$H_{3}(2)(f^{-1})=\frac{4}{15}c_{2}c_{4}-\frac{1}{4}c_{3}^{2}-\frac{28}{45}c_{2}^{3},$$
and further,
$$|H_{3}(2)(f^{-1})|\leq\frac{4}{15}|c_{2}||c_{4}|+\frac{1}{4}|c_{3}|^{2}+\frac{28}{45}|c_{2}|^{3}.$$
Applying Lemma \ref{lem-carl} (with $c_{1}=0$ ) we receive
$$|H_{3}(2)(f^{-1})|\leq\frac{4}{15}|c_{2}|(1-|c_{2}|^{2})+\frac{1}{4}(1-|c_{2}|^{2})^{2}+\frac{28}{45}|c_{2}|^{3}.$$
and, finally,
\be\label{eq-11}
|H_{3}(2)(f^{-1})|\leq \frac{1}{4}+\frac{4}{15}|c_{2}|-\frac{1}{2}|c_{2}|^{2}
+\frac{16}{45}|c_{2}|^{3}+\frac{1}{4}|c_{2}|^{4}=:\varphi_{1}(|c_{2}|),
\ee
where $0\leq|c_{2}|\leq1 $.
Since
\[
\begin{split}
\varphi_{1}'(|c_{2}|)&=\frac{4}{15}-|c_{2}|+\frac{16}{15}|c_{2}|^{2}+|c_{2}|^{3}\\
&=\frac{4}{15}(1-2|c_{2}|)^{2}+\frac{1}{15}|c_{2}|+|c_{2}|^{3}>0,
\end{split}
\]
we have $\varphi_{1}(|c_{2}|)\leq \varphi_{1}(1)=\frac{28}{45}$, and from \eqref{eq-11},
$$|H_{3}(2)(f^{-1})|\leq\frac{28}{45}=0.622\ldots .$$

This result is best possible as the function $f_{1}(z)=\ln\frac{1+z}{1-z}-z $ defined by $f_{1}'(z)=\frac{1+z^{2}}{1-z^{2}},$ shows.

\medskip

\item[(b)]   We apply the same method as in the previous case. Namely, from the definition of the class $\mathcal{C}$  we have
$$1+\frac{zf''(z)}{f'(z)} = \frac{1+\omega(z)}{1-\omega(z)},$$
where $\omega$ is a Schwartz function, and from here
\be\label{eq-12}
(z f'(z))'=\left[1+2\left(\omega(z)+\omega^{2}(z)+\cdots\right)\right]\cdot f'(z).
\ee
Using the notations \eqref{eq-1} and \eqref{eq-7}, and comparing the coefficients  in the relation \eqref{eq-12}, after some simple calculations, we obtain
\be\label{eq-13}
\left\{\begin{split}
a_{2}&=c_{1}, \\
a_{3}&=\frac{1}{3}\left(c_{2}+3c_{1}^{2}\right),\\
a_{4}&=\frac{1}{6}\left(c_{3}+5 c_{1}c_{2}+6 c_{1}^{3}\right)\\
a_{5}&=\frac{1}{30}\left(3c_{4}+14c_{1}c_{3}+43c_{1}^{2}c_{2}+30c_{1}^{4}+6c_{2}^{2}\right).
\end{split}
\right.
\ee
If $a_{2}=0$, then by \eqref{eq-13} we have $c_{1}=0$ , which implies
\be\label{eq-14}
a_3= \frac{1}{3}c_2,\quad a_4= \frac{1}{6}c_3,\quad a_5= \frac{1}{10}(c_4+2c_2^2).
\ee
Using \eqref{eq-5} and \eqref{eq-14} we obtain
$$H_{3}(2)(f^{-1})=\frac{1}{180}\left(6c_{2}c_{4}-5c_{3}^{2}-8c_{2}^{3}\right).$$
From the last relation we get
$$|H_{3}(2)(f^{-1})|\leq\frac{1}{180}\left(6|c_{2}||c_{4}|+5|c_{3}|^{2}+8|c_{2}|^{3}\right),$$
and further, after applying Lemma1 (with $c_{1}=0$ ),
$$|H_{3}(2)(f^{-1})|\leq\frac{1}{180}\left(6|c_{2}|(1-|c_{2}|^{2})+5(1-|c_{2}|^{2})^{2}+8|c_{2}|^{3}\right),$$
i.e.,
\be\label{eq-15}
|H_{3}(2)(f^{-1})|\leq \frac{1}{180}\left(5+6|c_{2}|-10|c_{2}|^{2}+2|c_{2}|^{3}+5|c_{2}|^{4}\right) =:\varphi_{2}(|c_{2}|),
\ee
where $0\leq|c_{2}|\leq1 $. Since
$$ \varphi_{2}'(|c_{2}|)=\frac{1}{90}\left(3-10|c_{2}|+3|c_{2}|^{2}+10|c_{2}|^{3}\right),$$
which, after considering this polynomial in the interval $[0,1]$, gives $\varphi_{1}(|c_{2}|)\leq \varphi_{2}(1)=\frac{2}{45}$, and further, from \eqref{eq-15},
$$|H_{3}(2)(f^{-1})|\leq\frac{2}{45}=0.044\ldots .$$
The function $f_{2}(z)=\operatorname{artanh}z$ satisfying $1+\frac{zf_2''(z)}{f'_2(z)}=\frac{1+z^2}{1-z^2}$ shows that the result is the best possible.

\medskip

\item[(c)] From the definition of the class $\mathcal{S}^{\star}$  we have that there exists a
Schwartz function $\omega$ such that
$$\frac{zf'(z)}{f(z)} = \frac{1+\omega(z)}{1-\omega(z)},$$
and from here
\be\label{eq-16}
z f'(z)=\left[1+2\left(\omega(z)+\omega^{2}(z)+\cdots\right)\right]\cdot f(z).
\ee
As in the two previous cases ((a) and (b)), by comparing the coefficients  in the relation \eqref{eq-16}, and some simple calculations, we have
\[
\left\{\begin{split}
a_{2}&=2c_{1} \\
a_{3}&=c_{2}+3c_{1}^{2}\\
a_{4}&=\frac{2}{3}\left(c_{3}+5 c_{1}c_{2}+6 c_{1}^{3}\right)\\
a_{5}&=\frac{1}{2}\left(c_{4}+\frac{14}{3}c_{1}c_{3}+\frac{43}{3}c_{1}^{2}c_{2}+10c_{1}^{4}+2c_{2}^{2}\right).
\end{split}
\right.
\]
For the case $a_{2}=0$ we have the next
\be\label{eq-18}
a_3= c_2,\quad a_4= \frac{2}{3}c_3,\quad a_5= \frac{1}{2}(c_4+2c_2^2).
\ee
So, from \eqref{eq-5} and \eqref{eq-18} we obtain
$$H_{3}(2)(f^{-1})=\frac{1}{18}\left(9c_{2}c_{4}-8c_{3}^{2}-36c_{2}^{3}\right),$$
and from here
$$|H_{3}(2)(f^{-1})|\leq\frac{1}{18}\left(9|c_{2}||c_{4}|+8|c_{3}|^{2}+36|c_{2}|^{3}\right).$$
Using estimates for $|c_{4}|$ and $|c_{3}|$ from Lemma1 (with $c_{1}=0$ ) from the last relation we receive
\be\label{eq-19}
|H_{3}(2)(f^{-1})|\leq \frac{1}{18}\left(8+9|c_{2}|-16|c_{2}|^{2}+27|c_{2}|^{3}+8|c_{2}|^{4}\right) =:\varphi_{3}(|c_{2}|),
\ee
where $0\leq|c_{2}|\leq1 $.
Since
\[
\begin{split}
\varphi_{3}'(|c_{2}|)&=\frac{1}{18}\left(9-32|c_{2}|+81|c_{2}|^{2}+32|c_{2}|^{3}\right)\\
&= \frac{1}{18}\left[9(1-3|c_{2}|)^{2}+22|c_{2}|+32|c_{2}|^{3}\right]>0,
\end{split}
\]
then $\varphi_{3}(|c_{2}|)\leq \varphi_{3}(1)=2$, and from \eqref{eq-19},
$$|H_{3}(2)(f^{-1})|\leq 2 .$$
The result is the best possible as the function $f_{3}(z)=\frac{z}{1-z^{2}}$ shows.

\medskip

\item[(d)] From the definition of the class $\mathcal{S}^{\star}_{s}$  we have that there exists a
Schwartz function $\omega$ such that
$$\frac{2zf'(z)}{f(z)-f(-z)} = \frac{1+\omega(z)}{1-\omega(z)},$$
and from here
\be\label{eq-20}
2zf'(z)=\left[1+2(\omega(z)+\omega^{2}(z)+\cdots)\right]\cdot[f(z)-f(-z)].
\ee
Similarly as in  previous cases, by comparing the coefficients  in the relation \eqref{eq-20}, after some simple calculations, we receive
\be\label{eq-21}
\left\{\begin{split}
a_{2}&=c_{1} \\
a_{3}&=c_{2}+c_{1}^{2}\\
a_{4}&=\frac{1}{2}\left(c_{3}+3c_{1}c_{2}+2c_{1}^{3}\right)\\
a_{5}&=\frac{1}{2}\left(c_{4}+2c_{1}c_{3}+5c_{1}^{2}c_{2}+2c_{1}^{4}+2c_{2}^{2}\right).
\end{split}
\right.
\ee
For $a_{2}=0$ ($c_1=0$), from \eqref{eq-21} we get
\[
a_3= c_2,\quad a_4= \frac{1}{2}c_3,\quad  a_5= \frac{1}{2}(c_4+2c_2^2),
\]
and using \eqref{eq-5},
$$H_{3}(2)(f^{-1})=\frac{1}{4}\left(2c_{2}c_{4}-c_{3}^{2}-8c_{2}^{3}\right),$$
and from here
$$|H_{3}(2)(f^{-1})|\leq\frac{1}{4}\left(2|c_{2}||c_{4}|+|c_{3}|^{2}+8|c_{2}|^{3}\right).$$
Using the estimates for $|c_{4}|$ and $|c_{3}|$ from Lemma \ref{lem-carl} (with $c_{1}=0$)  from the last relation we have
\be\label{eq-23}
|H_{3}(2)(f^{-1})|\leq \frac{1}{4}\left(1+2|c_{2}|-2|c_{2}|^{2}+6|c_{2}|^{3}+|c_{2}|^{4}\right) =:\varphi_{4}(|c_{2}|),
\ee
where $0\leq|c_{2}|\leq1 $.
Since
\[
\begin{split}
 \varphi_{4}'(|c_{2}|)&=\frac{1}{2}\left(1-2|c_{2}|+9|c_{2}|^{2}+2|c_{2}|^{3}\right)\\
 &=\frac{1}{2}\left[(1-|c_{2}|)^{2}+8|c_{2}|+2|c_{2}|^{3}\right]>0,
 \end{split}
 \]
then $\varphi_{4}$ is an increasing function and  $\varphi_{4}(|c_{2}|)\leq \varphi_{4}(1)=2$. So, from \eqref{eq-23},
$$|H_{3}(2)(f^{-1})|\leq 2 .$$
This result is the best possible as the function $f_{4}$ defined by
$$\frac{2zf_{4}'(z)}{f_{4}(z)-f_{4}(-z)} = \frac{1+z^{2}}{1-z^{2}}$$
shows.
\end{itemize}
\end{proof}

\medskip

\begin{rem}
From the relation \eqref{eq-6} we get the following.
\begin{itemize}
  \item[(a)] For $f\in\mathcal{R}$,
$$|H_{3}(2)(f^{-1})-H_{3}(2)(f)|=3|a_{3}|^{3}=3\left(\frac23|c_2|\right)^3\leq\frac{8}{9},$$
and the result is the best possible as the function $f_{1}$ shows (in this case $H_{3}(2)(f_{1})=\frac{4}{15}$ and $H_{3}(2)(f^{-1}_{1})=-\frac{28}{45}$).
  \item[(b)] For $f\in\mathcal{C}$,
$$|H_{3}(2)(f^{-1})-H_{3}(2)(f)|=3|a_{3}|^{3}=3\left(\frac{|c_2|}{3}\right)^3\leq\frac{1}{9},$$
and the result is the best possible as the function $f_{2}$ shows.
  \item[(c)] For $f\in\mathcal{S}^{\star}$,
$$|H_{3}(2)(f^{-1})-H_{3}(2)(f)|=3|a_{3}|^{3}=3|c_2|^3\leq3,$$
and the result is the best possible as the function $f_{3}$ shows.
  \item[(d)] For $f\in\mathcal{S}^{\star}_{s}$,
$$|H_{3}(2)(f^{-1})-|H_{3}(2)(f)|=3|a_{3}|^{3}=3|c_2|^3\leq 3,$$
and the result is the best possible for the function $f_{4}$.
\end{itemize}
\end{rem}

\medskip

For obtaining the corresponding result for the whole class $\mathcal{S}$  we will use method based on Grunsky coefficients.
In the proof we will use mainly the notations and results given in the book of N.A. Lebedev (\cite{Lebedev}).

Here are basic definitions and results.

Let $f \in \mathcal{S}$ and let
\[
\log\frac{f(t)-f(z)}{t-z}=\sum_{p,q=0}^{\infty}\omega_{p,q}t^{p}z^{q},
\]
where $\omega_{p,q}$ are the Grunsky's coefficients with property $\omega_{p,q}=\omega_{q,p}$.
For those coefficients we have the next Grunsky's inequality (\cite{duren,Lebedev}):
\be\label{eq-24}
\sum_{q=1}^{\infty}q \left|\sum_{p=1}^{\infty}\omega_{p,q}x_{p}\right|^{2}\leq \sum_{p=1}^{\infty}\frac{|x_{p}|^{2}}{p},
\ee
where $x_{p}$ are arbitrary complex numbers such that last series converges.

Further, it is well-known that if the function $f$ given by \eqref{eq-1}
belongs to $\mathcal{S}$, then also
\be\label{eq-25}
\tilde{f_{2}}(z)=\sqrt{f(z^{2})}=z +c_{3}z^3+c_{5}z^{5}+\cdots
\ee
belongs to the class $\mathcal{S}$. Then, for the function $\tilde{f_{2}}$ we have the appropriate Grunsky's
coefficients of the form $\omega_{2p-1,2q-1}^{(2)}$ and the inequality \eqref{eq-24} has the form:
\be\label{eq-26}
\sum_{q=1}^{\infty}(2q-1) \left|\sum_{p=1}^{\infty}\omega_{2p-1,2q-1}x_{2p-1}\right|^{2}\leq \sum_{p=1}^{\infty}\frac{|x_{2p-1}|^{2}}{2p-1}.
\ee

Here, and further in the paper we omit the upper index (2) in  $\omega_{2p-1,2q-1}^{(2)}$ if compared with Lebedev's notation.

\medskip

If in the inequality \eqref{eq-26} we put $x_{1}=1$ and $x_{2p-1}=0$ for $p=2,3,\ldots$, then we receive
\be\label{eq-27}
|\omega_{11} |^2 +3|\omega_{13}|^2 + 5|\omega_{15} |^2 +7|\omega_{17}|^2\leq 1.
\ee

As it has been shown in \cite[p.57]{Lebedev}, if $f$ is given by \eqref{eq-1} then the coefficients $a_{2}$, $ a_{3}$, $ a_{4}$ and $a_5$ are expressed by Grunsky's coefficients  $\omega_{2p-1,2q-1}$ of the function $\tilde{f}_{2}$ given by
\eqref{eq-25} in the following way:
\be\label{eq-28}
\begin{split}
a_{2}&=2\omega _{11},\\
a_{3}&=2\omega_{13}+3\omega_{11}^{2}, \\
a_{4}&=2\omega_{33}+8\omega_{11}\omega_{13}+\frac{10}{3}\omega_{11}^{3},\\
a_{5}&=2\omega_{35}+8\omega_{11}\omega_{33}+5\omega_{13}^{2}+18\omega_{11}^2\omega_{13}+\frac73\omega_{11}^4,\\
0&= 3\omega_{15}-3\omega_{11}\omega_{13}+\omega_{11}^3-3\omega_{33},\\
0&=\omega_{17}-\omega_{35}-\omega_{11}\omega_{33}-\omega_{13}^{2}+\frac{1}{3}\omega_{11}^{4}.
\end{split}
\ee

We note that in the cited book of Lebedev there exists a typing mistake for the coefficient $a_{5}$. Namely, instead of the term $5\omega_{13}^{2}$, there is $5\omega_{15}^{2}$.

\bthm\label{22-th-2}
Let $f\in\mathcal{S}$ is given by \eqref{eq-1} and let $a_{2}=0$. Then
 $$|H_{3}(2)(f^{-1})|\leq \frac{\sqrt{3}}{6\sqrt{7}}+2\sqrt{3}=3.57321\ldots.$$
\ethm

\begin{proof}
In the case when $a_{2}=0$, from \eqref{eq-28} we have $\omega _{11}=0$, and so
\be\label{eq-29}
a_{3}=2\omega_{13},\quad  a_{4}=2\omega_{33},\quad a_{5}=2\omega_{35}+5\omega_{13}^{2},\quad \omega_{33}=\omega_{15},\quad
\omega_{35}=\omega_{17}-\omega_{13}^{2}.
\ee
Using \eqref{eq-5} and \eqref{eq-29}, we have
$$H_{3}(2)(f^{-1})=4\omega_{13}\omega_{35}-14\omega_{13}^{3}-4\omega_{33}^{2},$$
and after applying the two last relations from \eqref{eq-29},
$$H_{3}(2)(f^{-1})=4\omega_{13}\omega_{17}-18\omega_{13}^{3}-4\omega_{15}^{2}.$$
From here we have
$$|H_{3}(2)(f^{-1})|\leq 4|\omega_{13}||\omega_{17}|+18|\omega_{13}|^{3}+4|\omega_{15}|^{2},$$
or finally, using  $|\omega_{17}|\leq \frac{1}{\sqrt{7}}\sqrt{1-3|\omega_{13}|^{2}-5|\omega_{15}|^{2}}$
(from \eqref{eq-27}) we get
\be\label{eq-30}
\begin{split}
|H_{3}(2)(f^{-1})|&\leq \frac{1}{\sqrt{7}}|\omega_{13}|\sqrt{1-3|\omega_{13}|^{2}-5|\omega_{15}|^{2}}+18|\omega_{13}|^{3}+4|\omega_{15}|^{2}\\
&=:\frac{1}{\sqrt{7}}\psi_{1}(|\omega_{13}|,|\omega_{15}|)+2\psi_{2}(|\omega_{13}|,|\omega_{15}|),
\end{split}
\ee
where
\[\psi_{1}(y,z)=y\sqrt{1-3y^{2}-5z^{2}},\qquad \psi_{2}(y,z)=9y^{3}+2z^{2},
\]
with $0\leq y=|\omega_{13}| \leq \frac{1}{\sqrt{3}}$ and
$0\leq z=|\omega_{15}| \leq\frac{1}{\sqrt{5}}\sqrt{1-3y^{2}}$
(where we used the inequality \eqref{eq-27}).
It is easy to verify that for these range of $y$ and $z$, $\psi_{1}(y,z)\leq\psi_1(1/\sqrt6,0)=\frac{\sqrt{3}}{6}$ and
$\psi_{2}(y,z)\leq\psi_2(1/\sqrt3,0)=\sqrt{3}$, so that from \eqref{eq-30}) we have
$$ |H_{3}(2)(f^{-1})|\leq \frac{\sqrt{3}}{6\sqrt{7}}+2\sqrt{3}=3.57321\ldots.$$
\end{proof}

\begin{rem}
From the relation \eqref{eq-6} we get for $f\in\mathcal{S}$:
$$|H_{3}(2)(f^{-1})-H_{3}(2)(f)|=3|a_{3}|^{3}= 3|2\omega_{13}|^{3}\leq 3\left(2\cdot \frac{1}{\sqrt3}\right)^3 =\frac{8}{\sqrt{3}}=4.6188\ldots. $$
\end{rem}

\end{document}